\documentclass{article}      % Specifies the document class \usepackage{epsfig}

\usepackage{epsfig}

\title{Q-Damped Oscillator and Degenerate Roots of Constant Coefficients q-Difference ODE}

\author{Sengul Nalci and Oktay K. Pashaev \\Department of Mathematics, Izmir Institute of Technology \\ Urla-Izmir, 35430, Turkey}

\begin{document}
\newcommand{\be}{\begin{equation}}
\newcommand{\ee}{\end{equation}}
\newcommand{\bea}{\begin{eqnarray}}
\newcommand{\eea}{\end{eqnarray}}
\newcommand{\disp}{\displaystyle}
\newcommand{\la}{\langle}
\newcommand{\ra}{\rangle}

\newtheorem{thm}{Theorem}[subsection]
\newtheorem{cor}[thm]{Corollary}
\newtheorem{lem}[thm]{Lemma}
\newtheorem{prop}[thm]{Proposition}
\newtheorem{defn}[thm]{Definition}
\newtheorem{rem}[thm]{Remark}
\newtheorem{prf}[thm]{Proof}

\maketitle

%%%%%%%%%%%%%%%%%%%%%%%%%%%%%%%%%%%%%%%%%%%%%%%%%%%%%%%%%%%%%% % You may repeat \author \address as often as necessary      % %%%%%%%%%%%%%%%%%%%%%%%%%%%%%%%%%%%%%%%%%%%%%%%%%%%%%%%%%%%%%%

\begin{abstract}
The classical model of q-damped oscillator is introduced and solved in terms of Jackson q-exponential function for three different cases, under-damped, over-damped and the critical one. It is shown that in all three cases solution is oscillating in time but is unbounded and non-periodic. By q-periodic function modulation, the self-similar microstructure of the solution for small time intervals is derived.
In the critical case with degenerate roots, the second linearly independent solution is obtained as a limiting case of two infinitesimally close roots. It appears as standard derivative of q-exponential and is rewritten
in terms of the q-logarithmic function. We extend our result by constructing n linearly independent set of solutions to a generic constant coefficient q-difference equation degree N with n degenerate roots.
\end{abstract}

\section{Introduction}

The damped harmonic oscillator as simplest classical model of motion with dissipation, corresponds to friction force proportional
to velocity of motion. It appears in many physical problems from quantum theory to inflating universe models.
The quantum damped oscillator as one of the simplest quantum system displaying the energy dissipation, has been studied to understand dissipation in quantum theory \cite{Eti and Pashaev}. Most popular models are, the Bateman-Feshbach-Tikochinsky oscillator as a closed system with two degrees of freedom and the Caldirola-Kanai oscillator as an open system with one degree of freedom and time dependent mass.
Due to complicated character of the friction force, several modifications of the damping term were proposed as fractional derivative, time-delay or finite difference derivative etc. The goal of the present paper is to study classical q-extended damped harmonic oscillator, where standard time derivatives are replaced by Jackson q-derivatives \cite{Kac et al.}. The q-extension of harmonic oscillator and its solution in form of basic trigonometric functions \cite{Jackson}
was considered before in \cite{Exton}. The q-extended heat equation and corresponding Burgers equation with q-shock solitons were studied in \cite{Nalci et al.}.
In the limit $q \rightarrow 1$, the q-deformed model reduces to the standard damped oscillator model. We construct solution in terms of Jackson q-exponential function and find q-periodic modulation of the solution with self-similar properties. Special attention is paid for degenerate roots case.
And results are generalized for arbitrary order constant coefficient q-ordinary difference equation.
 This gives background for further possible quantization of corresponding model.

\section{Damped Oscillator}

In reality a spring never oscillates forever, since frictional forces will diminish the amplitude of oscillation until the rest.  In many situations the frictional force is proportional to the velocity of the mass as follows
$f_r= -\gamma v,$ where $\gamma >0$ is the damping constant. Therefore, by adding this frictional force we have the following equation for a spring
\be m \frac{d^2 x}{d t^2}+ \gamma \frac{d x}{d t}+ kx=0. \label{dampedos} \ee
Solution of this equation in the form $x(t)= e^{\lambda t}$ leads to the characteristic equation
\be m \lambda^2 +\gamma \lambda +k =0 \ee with two roots
$$\lambda_1= \frac{-\gamma + \sqrt{\gamma^2-4 m k}}{2 m}, \,\,\,\,\,\,\,\,\,\, \lambda_2= \frac{-\gamma - \sqrt{\gamma^2-4 m k}}{2 m}.$$
Then according to value of damping constant we have three cases : \\
\textbf{i\,- Under-damping Case: } When $\gamma^2< 4 m k ,$ which means that friction is sufficiently weak, we have two complex conjugate roots
\be \lambda_{1,2}= -\frac{\gamma}{2 m} \pm i \omega ,\ee
where $\omega \equiv \sqrt{\frac{k}{m}-\frac{\gamma^2}{4 m^2}}.$ Then the general solution of (\ref{dampedos}) is
\be x(t)= e^{-\frac{\gamma}{2 m}t}\left(A \cos \omega t+ B\sin \omega t\right).\ee
If $\gamma=0,$ there is no decay and the spring oscillates forever.
If $\gamma$ is big, the amplitude of oscillations decays very fast (the exponential decay).\\
\textbf{ii\,- Over-damping Case: } When $\gamma^2> 4 m k ,$ which means that friction is sufficiently strong,  both roots are real,\, this why the solution decays  exponentially
\be x(t)= A e^{\frac{-\gamma + \sqrt{\gamma^2-4m k}}{2 m}t}+ B e^{\frac{-\gamma - \sqrt{\gamma^2-4m k}}{2 m}t}.\ee
This case is called as over-damping because there is no any oscillation.\\
\textbf{iii\,- Critical Case: } For $\gamma^2= 4 m k,$ we have two degenerate roots
$$\lambda_1=\lambda_2= -\frac{\gamma}{2 m},$$ then  the general solution is
\be x(t)= A e^{-\frac{\gamma}{2 m}t}+ B t e^{-\frac{\gamma}{2 m}t}.\ee
\section{q-Harmonic Oscillator}
Here we introduce the $q$-Harmonic oscillator.
Equation of  $q$-deformed classical harmonic oscillator is
\be D_q^2 x(t)+ \omega^2 x(t)=0,\ee
where the q-derivative is definite as \cite{Kac et al.},
\be  D_q x(t) = \frac{x(qt) - x(t)}{(q-1)t}.\ee
Using the power series method (or the $q$-exponential form $x(t)=e_q(\lambda t)$), we find the general solution of $q$-Harmonic Oscillator in the following form
\cite{Exton}
\be x(t)= A(t) \cos_q \omega t+ B(t) \sin_q \omega t,\ee where $$D_q A(t)= D_q B(t)=0,$$ means $A(t),B(t)$ in general are $q$-periodic functions, and particularly could be arbitrary constants. Here the Jackson q-exponential function is definite as
\be e_q (t)   \sum^\infty_{n=1} \frac{t^n}{[n]_q !},\ee
and
\be e_q (it) = \cos_q t + i \sin_q t, \ee
where $[n]_q = 1+ q + ... + q^{n-1}$. For $q> 1$ this function is entire analytic function, so we restrict consideration by this case only.

\begin{figure}[h]
\begin{center}
{\includegraphics[width=3.5in]{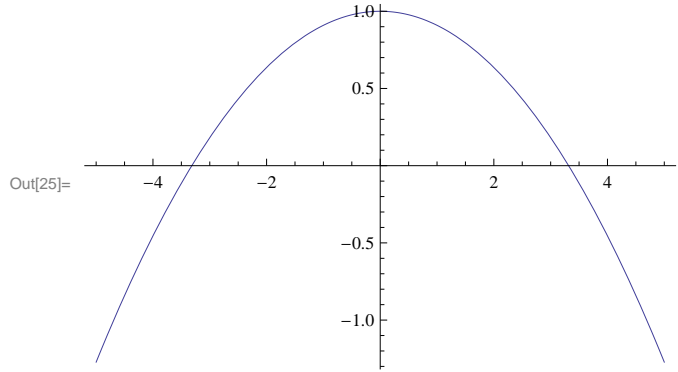}\\}
\caption{ q-Harmonic oscillator solution $\cos_q t$}\label{cosq1}
\end{center}
\end{figure}

\begin{figure}[ht]
\begin{center}
{\includegraphics[width=3.5in]{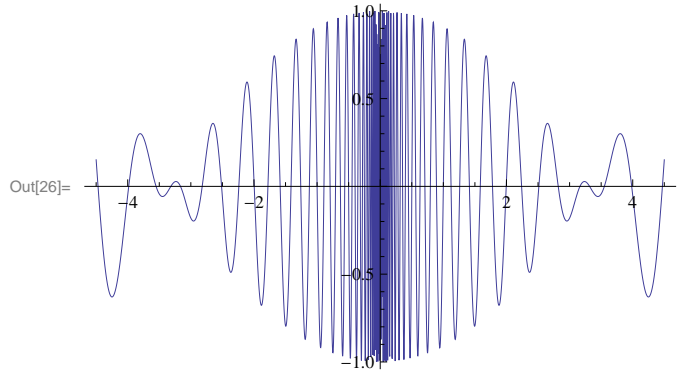}\\}
\caption{ q-Harmonic oscillator solution $\sin \left( \frac{2 \pi }{\ln q} \ln t\right)\cos_q t$}\label{cosq1periodic}
\end{center}
\end{figure}

In Figure \ref{cosq1} we plot particular $\cos_q t$ solution of $q$-deformed classical harmonic oscillator.
In contrast to standard $\sin t$ and $\cos t$ functions, $\sin_q t$ and $\cos_q t$ functions \cite{Jackson},
are not bounded and also have no periodicity.
In Figure \ref{cosq1periodic} we plot modulation of the same solution with q-periodic function $A(t)=\sin \left( \frac{2 \pi }{\ln q} \ln t\right)\cos_q t,$ which gives micro oscillations to the solution.

\newpage

\section{q-Damped Harmonic Oscillator}
We define equation for $q$-analogue of damped oscillator in the form
\be D_q^2 x(t)+ \Gamma D_q x(t)+ \omega^2 x(t)=0, \label{qdamped}\ee
where $$\omega \equiv \sqrt{\frac{k}{m}},\,\,\,\,\,\,\,\,\, \Gamma \equiv \frac{\gamma}{m}.$$
 By substituting $x(t)= e_q(\lambda t)$  into  equation (\ref{qdamped}), we obtain
\be e_q(\lambda t)\left[\lambda^2 +\Gamma \lambda +\omega^2 \right]=0.\ee
For $q>1,$ $e_q(\lambda t)$ is an entire function defined for any $t$, this why it has an infinite set of zeros (no poles). Then,\,we can choose
$$\lambda^2 +\Gamma \lambda +\omega^2 =0.$$ The roots of this characteristic equation are $$\lambda_{1,2}=-\frac{\Gamma}{2} \pm \sqrt{\frac{\Gamma^2}{4}-\omega^2}.$$

\subsection{Under-Damping Case}
For $\Gamma^2 < 4 \omega^2,$ we have two complex conjugate roots
$$\lambda_1 = -\frac{\Gamma}{2}+ i \Omega ,\,\,\,\,\,\,\,\, \lambda_2 = -\frac{\Gamma}{2}- i \Omega,$$
where $$\Omega \equiv \sqrt{\omega ^2-\frac{\Gamma^2}{4}}.$$
Then the general solution of equation (\ref{qdamped})\,is
\be x(t)= A e_q \left[\left(-\frac{\Gamma}{2}+ i \Omega \right)t\right]+B e_q \left[\left(-\frac{\Gamma}{2}-i \Omega \right)t\right]  \ee

\begin{figure}[ht]
\begin{center}
{\includegraphics[width=3.5in]{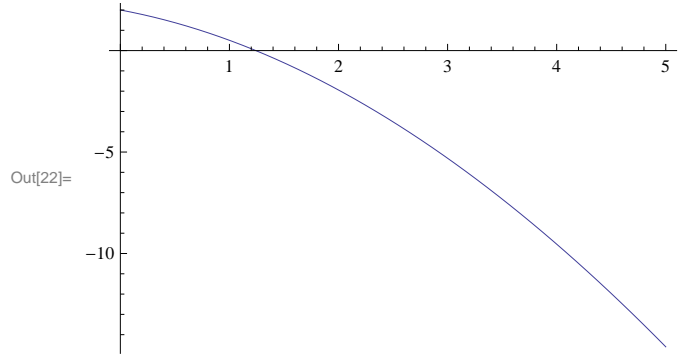}\\}
\caption{Under-damping case $A=B=1$}\label{qunderdamping}
\end{center}
\end{figure}

\begin{figure}[h]
\begin{center}
{\includegraphics[width=3.5in]{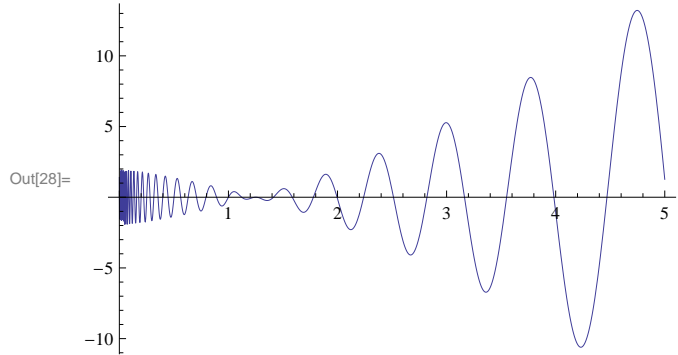}\\}
\caption{ Under-damping case with $q$-periodic function}\label{qunderdampingperiodic}
\end{center}
\end{figure}

In Figure \ref{qunderdamping} and Figure \ref{qunderdampingperiodic} we plot particular solutions with constant $(A=B=1)$  and with $q$-Periodic modulation, respectively.

\subsection{Over-Damping Case}
For $\Gamma^2 > 4 \omega^2,$ we have two distinct real roots $\lambda_{1,2}$  and solution is
\be x(t)= A(t) e_q\left[\left(-\frac{\Gamma}{2} + \sqrt{\frac{\Gamma^2}{4}-\omega^2}\right)t\right]+B(t)e_q\left[\left(-\frac{\Gamma}{2}- \sqrt{\frac{\Gamma^2}{4}-\omega^2}\right)t\right] ,\ee
where $A(t),B(t)$ are $q$-periodic functions (or could be arbitrary constants).

\begin{figure}[h]
\begin{center}
{\includegraphics[width=3.5in]{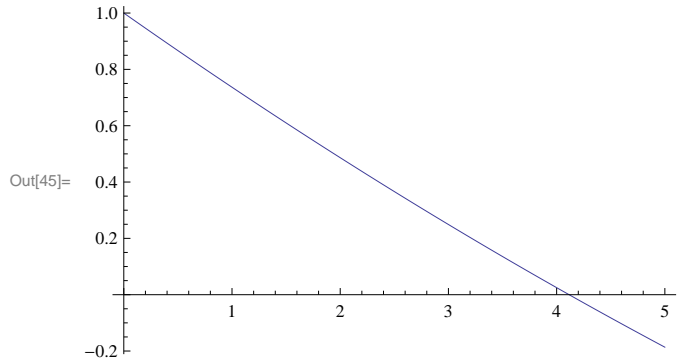}\\}
\caption{Over-damping case $A=B=1$}\label{qoverdamping}
\end{center}
\end{figure}

\begin{figure}[ht]
\begin{center}
{\includegraphics[width=3.5in]{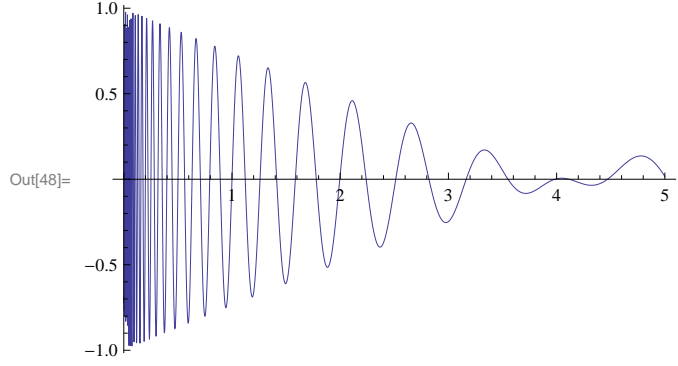}\\}
\caption{Over-damping case with $q$-periodic function}\label{qoverdampingperiodic}
\end{center}
\end{figure}

In Figure \ref{qoverdamping} and Figure \ref{qoverdampingperiodic} we plot particular solutions with constant $(A=B=1)$  and with $q$-Periodic modulation, respectively.

\subsection{Critical Case}
For $\Gamma^2 = 4 \omega^2,$ we have degenerate roots $\lambda_{1,2}=-\frac{\Gamma}{2}.$
The first obvious solution is $e_q(-\omega t).$ However if we try the second linearly independent solution in the usual form $t e_q(-\omega t),$ it doesn't work. This why we follow the next method:

 We suppose that the system is very close to the critical case so that $\frac{\Gamma}{2}= \omega + \epsilon,$
where $\epsilon \ll 1.$
Then the roots of characteristic equation are
\be \lambda_1= -\omega + \sqrt{2 \omega \epsilon}, \,\,\,\,\,\,\,\,\,\lambda_2= -\omega - \sqrt{2 \omega \epsilon}, \ee
and the solution is
\be x(t)= A e_q\left( (-\omega + \sqrt{2 \omega \epsilon} )t\right)+ B e_q\left( (-\omega - \sqrt{2 \omega \epsilon} )t\right)\ee
Expanding this solution in terms of $\epsilon,$
\bea x(t)&=& A \sum_{n=0}^\infty \left(\frac{(-\omega)^n + n (\sqrt{2 \omega \epsilon} (-\omega)^{n-1}+...)}{[n]!}\right)t^n \nonumber \\
& +& B \sum_{n=0}^\infty \left(\frac{(-\omega)^n - n (\sqrt{2 \omega \epsilon} (-\omega)^{n-1}+...)}{[n]!}\right)t^n  \nonumber \\
&=& (A+B) \sum_{n=0}^\infty \frac{(-\omega)^n}{[n]!} t^n + (A-B) \sqrt{2 \omega \epsilon}\sum_{n=1}^\infty \frac{n}{[n]!} (-\omega)^{n-1} t^n +... \nonumber \\
&=& (A+B) x_1(t)+ (B-A)\sqrt{\frac{2 \epsilon}{\omega}} \,\,x_2(t)+..., \eea
in zero approximation we get the first solution \be x_1(t)= e_q(-\omega t).\ee
In the linear approximation we obtain the second solution in the form
\be x_2(t)= t \frac{d}{d t} e_q(-\omega t).\ee

In order to prove that solutions $x_1(t)= e_q(-\omega t)$ and $x_2(t)= t \frac{d}{dt} e_q(-\omega t)$ are linearly independent, we check the $q$- Wronskian :
\[ W_q= \left|
\begin{array}{cc}
   e_q(-\omega t) &   t \frac{d}{dt} e_q(-\omega t)  \\
  D_q(e_q(-\omega t))     & D_q\left(t \frac{d}{dt} e_q(-\omega t)\right)
\end{array} \right| \],
or
\be W_q= -\omega e_q(-\omega t) \left( e_q(-\omega t)-t \frac{d}{d t}e_q(-\omega t)\right).\ee
Here we show that  the term in parenthesis is not identically zero.
For $q>1,$ by using the infinite product representation of $e_q(x)$ (\ref{infprod}), we get
\be e_q(-\omega t)= \prod_{n=0}^\infty \left(1-\left(1-\frac{1}{q}\right) \frac{1}{q^n} \omega t \right),\ee
$$t \frac{d}{dt} \ln e_q(-\omega t)= \sum_{n=0}^\infty \frac{-w \left(1-\frac{1}{q}\right)\frac{1}{q^n}t}{1-\left(1-\frac{1}{q}\right) \frac{1}{q^n} \omega t},$$
or
$$t \frac{d}{d t} e_q(-\omega t)= A e_q(-\omega t),$$ where $$A \equiv \sum_{n=0}^\infty \frac{-\omega \left(1-\frac{1}{q}\right)\frac{1}{q^n}t}{1-\left(1-\frac{1}{q}\right) \frac{1}{q^n} \omega t}.$$
Expanding the denominator, we have
\bea A&=& \sum_{n=0}^\infty \left( -\left(1-\frac{1}{q}\right)\frac{1}{q^n} \omega t\right) \sum_{l=0}^\infty \left(1-\frac{1}{q}\right)^l \frac{1}{q^{nl}} (\omega t)^l \nonumber \\
&=& - \sum_{l=0}^\infty \left( 1-\frac{1}{q}\right)^{l+1} (\omega t)^{l+1} \sum_{n=0}^\infty \frac{1}{q^{n(l+1)}} \nonumber \\
&=& -\sum_{l=1}^\infty \frac{\left( \left(1-\frac{1}{q}\right) \omega t\right)^l}{[l]} \frac{1}{1-q}
\eea
where $|t|< \frac{q}{w}.$
We know that
$$\ln(1-x)= -x-\frac{x^2}{2}-\frac{x^3}{3}-...= -\sum_{l=1}^\infty \frac{x^l}{l}$$ and the $q$-analogue of this expression is given as \cite{Pashaev et al.}
\be {\rm Ln}_q (1-x)=- \sum_{l=1}^\infty \frac{x^l}{[l]_q}. \ee Then we rewrite
\be A= -\frac{1}{1-q} Ln_q \left( 1-\left(1-\frac{1}{q}\right) \omega t\right). \ee

The second solution $x_2(t)$ can also be rewritten in terms of $q$-logarithmic function
\be x_2(t)= t \frac{d}{d t} e_q(-\omega t)= \frac{1}{q-1} Ln_q \left(1-\left(1-\frac{1}{q}\right) \omega t\right) e_q(-\omega t),\ee where $|t|< \frac{q}{\omega}.$
Finally, the $q$-Wronskian is not vanish
$$ W_q= -\omega \left(e_q (-\omega t)\right)^2 \left( 1- \frac{1}{q-1} Ln_q \left(1-\left(1-\frac{1}{q}\right) \omega t\right)\right)\neq 0 ,$$
since the term $$1- \frac{1}{q-1} Ln_q \left(1-\left(1-\frac{1}{q}\right) \omega t\right)$$ couldn't be identically zero.

We can also rewrite this in terms of $q$-logarithm, which instead of linear in $t$ term for $q=1$ case, now includes infinite set of arbitrary powers of $t,$
\bea x_2(t)&=& \frac{1}{1-q} \sum_{l=1}^\infty \frac{\left((1-\frac{1}{q}) \omega t\right)^l}{[l]} e_q(-\omega t) \nonumber \\
&=& -\frac{1}{1-q} Ln_q \left( 1-\left(1-\frac{1}{q}\right) \omega t\right) e_q(-\omega t).\eea
It is easy to check that for $q\rightarrow 1$ our solution reduces to the standard second solution $t e^{-\omega t}.$

Combining the above results we find the general solution in the degenerate case as
\be x(t)= A e_q(-\omega t)+ B\, t \frac{d}{d t}e_q(-\omega t).\ee

\begin{figure}[h]
\begin{center}
{\includegraphics[width=3.5in]{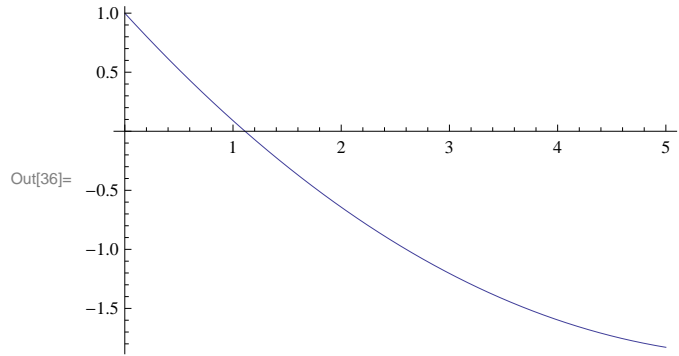}\\}
\caption{Critical case}\label{critical}
\end{center}
\end{figure}

\begin{figure}[ht]
\begin{center}
{\includegraphics[width=3.5in]{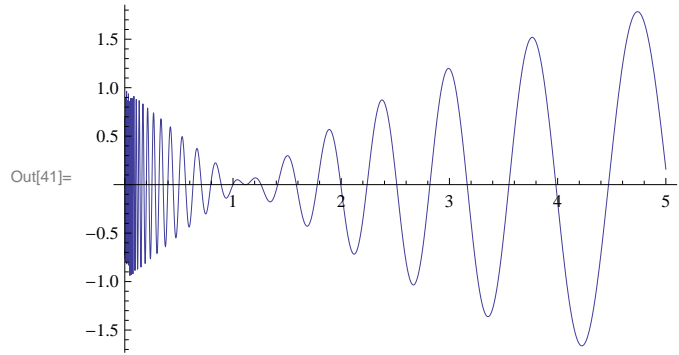}\\}
\caption{Critical case with periodic function}\label{criticalperiodic}
\end{center}
\end{figure}

In Figure \ref{critical} and Figure \ref{criticalperiodic} we plot particular solutions with constant $(A=B=1)$  and with $q$-Periodic modulation, respectively.

\begin{figure}[h]
\begin{center}
{\includegraphics[width=3.5in]{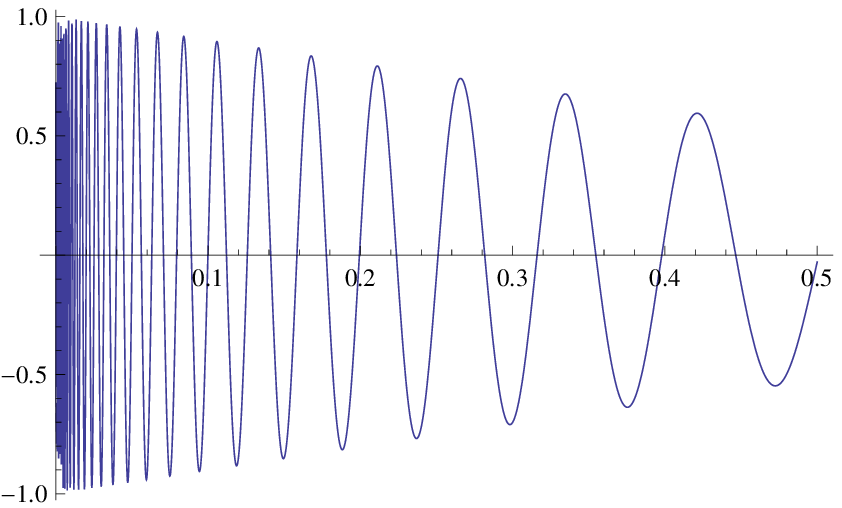}\\}
\caption{Self-similar micro structure at scale 0.5}\label{qperiodic1}
\end{center}
\end{figure}

\begin{figure}[ht]
\begin{center}
{\includegraphics[width=3.5in]{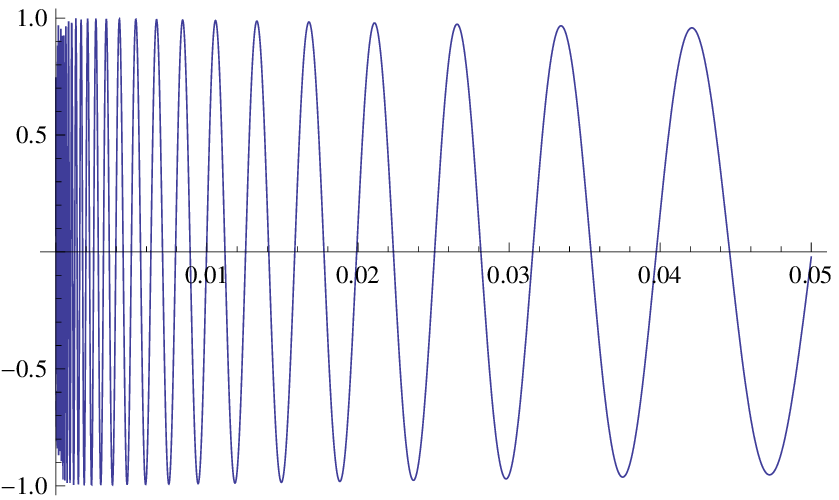}\\}
\caption{Self-similar micro structure at scale 0.05}\label{qperiodic2}
\end{center}
\end{figure}

In Figures \ref{qperiodic1} and \ref{qperiodic2} we plot solution with $q$-Periodic function modulation at different small scales.
Comparing these figures we find very close similarity, this why $q$-periodic function modulation leads to the self-similarity property of the solution.

\section{Degenerate Roots for Equation Degree N}

$q$-Damped oscillator considered in first section is an example of constant coefficient q-difference equation of degree two. The problem of radiation damping leads to a constant coefficient equation of degree three \cite{Eti and Pashaev}. Here we consider generic constant coefficient q-difference equation of degree N.
Then the result for degenerate roots obtained in previous section can be generalized to this equation of an arbitrary order.

The constant coefficients $q$-difference equation of order $N$ is
\be \sum_{k=0}^N a_k D^k x(t)=0, \label{constantcoefeqn}\ee where $a_k$ are constants (or $q$-periodic functions). By substitution \be x=e_q(\lambda t) \label{sol}\ee we get the characteristic equation  $$\sum_{k=1}^N a_k \lambda^k=0.$$  It has N roots.  Suppose $(\lambda_1, \lambda_2,...,\lambda_N)$ are distinct numbers.
Then, the general solution of (\ref{constantcoefeqn}) is found in the form
\be x(t)= \sum_{k=1}^N c_k e_q(\lambda_k t).\ee

In case, when we have $n$-degenerate roots $$(D+\omega)^n x=0,$$ by substituting (\ref{sol}), characteristic equation is found as $$(\lambda+\omega)^n =0. $$
Then the linearly independent solutions for these degenerate roots we can obtain in the following form :
$$x_1(t)= e_q(-\omega t),$$
$$x_2(t)= t \frac{d}{d t} e_q(-\omega t)= \frac{1}{q-1} {\rm Ln_q}\left(1-(1-\frac{1}{q}) \omega t \right) e_q(-\omega t), $$
$$x_3(t)=  \left(t \frac{d}{dt}\right)^2 e_q(-\omega t)$$ or  by using the commutation relation $\left[t, \frac{d}{d t}\right]=-1,$ up to linearly dependent solution, it can be written as
$$x_3(t)=t^2 \frac{d^2}{dt^2} e_q(-\omega t)$$
$$...$$
$$x_n(t)= t^{n-1} \frac{d^{n-1}}{dt^{n-1}} e_q(-\omega t).$$

For construction solutions with degenerate roots we need following propositions :

\begin{prop}
We have following commutation relation \be \left[ t \frac{d}{d t}, D \right]= -D \ee
implies $$t \frac{d}{d t}D= D\left( t \frac{d}{d t}-1\right)$$
\end{prop}
\begin{prf}
By definition of $D_q$ operator the commutation relation can be found as follows
\bea \left[ t \frac{d}{d t}, D \right]f &=& t \frac{d}{d t} D f- D t \frac{d}{d t} f  \nonumber \\
&=& t \frac{d}{d t} \left( \frac{f(q t)-f(t)}{(q-1)t}\right)- D \left( t \frac{d f}{dt}\right) \nonumber \\
&=& t \left( \frac{(q-1)t\left(q f'(q t)-f'(t)\right)-(q-1)\left(f(q t)-f(t) \right)}{(q-1)^2 t^2}\right)- \frac{q t \frac{d f(q t)}{d(q t)}-t f'(t)}{(q-1)t}, \nonumber \\
&=& -D f \eea which implies $$\left[ t \frac{d}{d t}, D \right]= -D .$$
\end{prf}

\begin{prop}
\be t \frac{d}{d t} D^n= D^n \left(t \frac{d}{d t}-n\right)\ee
\end{prop}

\begin{prf}
By using mathematical induction: \\
For $n=1,$ from the above commutation relation it is easy to see. \\
Suppose it is true for $n$:  $t \frac{d}{d t} D^n= D^n \left(t \frac{d}{d t}-n\right)$ \\
And we should show that it is true for $n+1$ \\
\bea t \frac{d}{d t} D^{n+1}&=&t \frac{d}{d t} D^n D = D^n \left( t \frac{d}{d t}-n\right) D = D^n t \frac{d}{d t} D -n D^{n+1} \nonumber \\
&=& D^n D \left(t \frac{d}{d t}-1 \right)  -n D^{n+1} = D^{n+1} \left(t \frac{d}{d t}-(n+1)\right) \eea
\end{prf}

\begin{prop}
We have more general relation in the following form
\be t \frac{d}{d t} (\omega+ D)^n= (\omega+D)^n t \frac{d}{d t}- n(\omega+ D)^{n-1} D \label{operatoridentity}.\ee
\end{prop}

\begin{prf}
By using mathematical induction : \\
For $n=1,$
\bea t \frac{d}{d t} (\omega+ D) &=& t \frac{d}{d t} \omega +t \frac{d}{d t} D = t \frac{d}{d t} \omega + D \left(t \frac{d}{d t}-1 \right) \nonumber \\
&=& (\omega+D) t \frac{d}{d t}- D  \nonumber \eea
Suppose this relation is true for $n$:
$$ t \frac{d}{d t} (\omega+ D)^n= (\omega+D)^n t \frac{d}{d t}- n(\omega+ D)^{n-1} D .$$
Now we prove that it is true for $n+1$
\bea  t \frac{d}{d t} (\omega+ D)^{n+1}&=& t \frac{d}{d t} (\omega+ D)^n (\omega+ D) \nonumber \\
&=& \left( (\omega+ D)^n t \frac{d}{d t}- n(\omega+ D)^{n-1} D \right)  (\omega+ D) \nonumber \\
&=& (\omega+ D)^n \left(  (\omega+ D)t \frac{d}{d t}-D \right)-n(\omega+ D)^n D \nonumber \\
&=& (\omega+ D)^{n+1} t \frac{d}{d t}-(n+1)(\omega+ D)^n D  \nonumber  \eea
\end{prf}

Using the operator identity (\ref{operatoridentity})
we can show that if $x_0$ is solution of $$(D+\omega) x_0=0 \Rightarrow (D+\omega)^2 x_0=0 \Rightarrow ...\Rightarrow (D+\omega)^n x_0=0.$$

\bea & & t \frac{d}{d t} (D+\omega)^n x_0=0 \nonumber \\
& & (D+\omega)^n t \frac{d}{d t} x_0- n D (\omega+ D)^{n-1} x_0 =0  \nonumber \eea
$$(\omega+ D)^{n-1} x_0 =0 \Rightarrow (\omega+ D)^n x_1=0,$$ where $$x_1\equiv  t \frac{d}{d t} x_0.$$

Then, $$x_1= t \frac{d}{d t} x_0$$ is solution of $$(D+\omega)^2 x_1=0,$$ and as follows $$(D+\omega)^n x_1=0,$$ e.t.c. And then,
$$x_{n-1}= t \frac{d^{n-1}}{d t^{n-1}} x_0$$ is solution of $$(D+\omega)^n x_{n-1}=0 .$$ It provides us with $n$ linearly independent solutions $x_0,x_1,...,x_{n-1}$ of $N$-degree equation with $n$-degenerate roots $$(D+ \omega)^n x=0.$$

\section{Conclusions}

In conclusion we like to mention relation of our q-damped oscillator problem with nonlinear q-difference problem.
By following substitution
\be y(t) = \frac{D_q x(t)}{x(t)} \ee
equation (\ref{qdamped}) leads to nonlinear the q-Riccati equation
\be D_q y(t) + y(qt) y(t) + \Gamma y(t) + \omega^2 = 0,  \label{riccati}\ee
so that every solution of the first one produces a solution of the second one. This means that (\ref{qdamped}) gives linearization
of the nonlinear q-difference equation  (\ref{riccati}). Similar situation we encounter in the q-Burgers equation
linearized by q-Cole-Hopf transformation in therms of the q-heat equation \cite{Nalci et al.}.

Finally some comments about degenerate limit of our system. If in addition to our q-damped oscillator model (\ref{qdamped})
we consider standard oscillator with q-derivative friction 
\be m\frac{d^2}{dt^2} x(t)+ \gamma D_q x(t)+ k x(t)=0, \label{sqdamped}\ee
then in the case of strong damping $\gamma >> m$, or the small mass, both equations reduce to the first order q-difference equation
\be \gamma D_q x(t)+ k x(t)=0. \label{dqdamped}\ee
Under the reduction, the two-dimensional phase space for the second order system turns into a one-dimensional one for the first  order system.
The first order system is called the "degenerate system" \cite{Andronov}. An arbitrary initial value problem in general does not apply
to the degenerate system (\ref{dqdamped})\cite{Vitiello}. For the second order system at $t=0$ we can attach an arbitrary value for coordinate $x$
and related velocity $\dot x$ or $D_q x$. However, we can describe the same physical system by (\ref{dqdamped}) only after some time interval and,  
moreover, in such a case  $\dot x$ or $D_q x$ cannot be arbitrary since $D_q x$ is completely determined by given coordinate $x$, according to (\ref{dqdamped}).
When the mass $m$ is going to zero, transition from a state incompatible with (\ref{dqdamped}) to a compatible one is very fast. Acceleration at the initial time is 
very high (the related velocity is changing very fast). The transition to the massless limit can be well approximated by the discontinuous jumping 
condition like in \cite{Andronov}: the energy of the system cannot be changed by a jump. The jumping condition implies that under the jumping the coordinates 
of the system remain invariant and only the velocities can be changed. Details of this study we present in our future work.

 \section*{Acknowledgments}
This work was support by
TUBITAK (The Scientific and Technological Research Council of Turkey), TBAG Project 110T679 and Izmir Institute of
Technology .

\end{document}